\documentclass{llncs}

\usepackage{bm}

\usepackage{amsmath,amsfonts,amssymb} 
\usepackage{graphicx} 

\def\x#1{} 

\newcommand{\one}{\underline{1}}

{\bfseries}{\itshape}
\newtheorem{lem}{Lemma}{\bfseries}{\itshape}
\newtheorem{cor}{Corollary}{\bfseries}{\itshape}
\newtheorem{pro}{Proposition}{\bfseries}{\itshape} 
{\bfseries}{\upshape}

\begin{document}

\title{Analysis of Relaxation Time in Random Walk with Jumps}

\author{Konstantin Avrachenkov$^1$ and Ilya Bogdanov$^2$}

\institute{Inria Sophia Antipolis, France\\ {\tt k.avrachenkov@inria.fr}
\and Higher School of Economics, Russia\\ {\tt ilya160897@gmail.com}}

\maketitle

\begin{abstract}
We study the relaxation time in the random walk with jumps. The random walk with
jumps combines random walk based sampling with uniform node sampling and improves
the performance of network analysis and learning tasks. We derive various conditions under
which the relaxation time decreases with the introduction of jumps.

\bigskip

\end{abstract}

\section{Introduction}

In the present work we study the relaxation time or equivalently the spectral gap
of the random walk with jumps on a general weighted undirected graph.
The random walk with jumps can be viewed as a random
walk on a combination of the original graph and the complete graph weighted by
a scaled parameter. This parameter determines the rate of jumps from the current
node to an arbitrary node. The random walk with jumps has similarities with
PageRank \cite{PR}. In fact, it coincides with PageRank
on the regular graphs but differs on the irregular graphs. In the case of the
random walk with jumps, the jump probability depends on the node degree.
The higher the degree of the current node is, the less likely the random walk
will jump out of the node to an arbitrary node. The random walk with jumps can also
be viewed as a particular case of the generalisation of PageRank with node-dependent
restart \cite{AHS14}.

The random walk with jumps has been introduced in \cite{ART10} to improve
the random walk based network sampling or respondent driven sampling \cite{VH08}
by combining the standard random walk based sampling with uniform node sampling.
A big advantage of the random walk with jumps in comparison with PageRank
is that, as opposite to PageRank, the random walk with jumps on an undirected graph
is a time-reversible process and its stationary distribution is available in a simple, explicit form.
In particular, this allows us to unbias efficiently the random walk, which is
node degree biased on irregular graphs. This comes with a price. The price is the
difficulty to control the relaxation time. In the case of PageRank the relaxation time
is bounded by the reciprocal of the restart probability \cite{HK03,LM06}.
In the case of the random walk with jumps, there is no simple connection with the jump parameter.
In \cite{ART10} the authors have shown that under a natural condition on the
clustering structure of the graph, the relaxation time decreases with the increase
of the jump parameter.

Let us mention a few more applications of the random walk with jumps beyond
network sampling \cite{ART10,Metal17,RT10}. The random walk with jumps has been used
in the context of graph-based semi-supervised learning \cite{ACM17,KL02}. The random
walk with jumps has also been used as a main building block in the quick algorithm
for finding largest degree nodes in a graph \cite{ALST14}. A continuous-time
version of the random walk with jumps has an application in epidemiology \cite{JT18}.

All this motivates us to take another look at the relaxation time of the random
walk with jumps. In particular, we are now able to give a necessary and sufficient
condition for the improvement of the relaxation time on weighted graphs. We give
an example showing that there are weighted graphs where introducing jumps increases
the relaxation time. The necessary and sufficient conditions are not easy to interpret.
Therefore we derive a series of simpler sufficient conditions. One new sufficient condition,
similar in spirit to the condition in \cite{ART10}, indicates that on graphs with clusters,
the relaxation time improves with the introduction of jumps.
The other new sufficient conditions require the spectral gap
of the original graph to be smaller than the reciprocal of the squared coefficient
of variation of the nodes' degrees, thus establishing a connection between the measure
of graph irregularity and the relaxation time. We expect that the derived conditions
are satisfied in most complex networks (either due to clustering structure or
due to small spectral gap). Thus, the present study confirms that it is safe and
in most cases beneficial to use the random walk with jumps for complex network analysis.

The structure of the paper is as follows: in the next section we define the random walk
with jumps and provide necessary background material. In Section~\ref{sec:dobr} we discuss
the application of Dobrushin coefficient for large jump rates. Then, in Section~\ref{sec:NDcond}
we discuss necessary and sufficient conditions for the reduction of the relaxation time when
the jump rate is small. In Section~\ref{sec:nodedegcond} we provide a series of sufficient conditions,
which are easier to interpret and to verify. In particular, we provide a sufficient condition in terms
of the coefficient of variation of nodes' degrees.
In Section~\ref{sec:num} we give interesting numerical illustrations.
We conclude the paper with a conjecture in Section~\ref{sec:conc}.

\section{Definitions and preliminaries}
\label{sec:def}

Most of the analysis in the present article is for the case of a general weighted
undirected graph $G$ with vertex set $V(G)$, $|V(G)|=n$, and edge set $E(G)$, $|E|=m$,
defined by the \textit{weighted adjacency matrix\/} $A=(a_{ij})$ with elements
$$
a_{ij}=\left\{\begin{array}{ll}
\mbox{weight of edge $(i,j)$}, & \mbox{if} \ i \sim j,\\
0, & \mbox{otherwise}.
\end{array}\right.
$$
Unless stated otherwise, we assume that $G$ is connected.

Denote by $\one$ the column vector of ones of appropriate dimension.
Then, $d=A\one$ is the vector of weighted degrees of vertices and $D=Diag(d)$
is the diagonal matrix with vector $d$ on the main diagonal.

The Standard Random Walk (SRW) is a discrete-time Markov chain $\{X_t, t=0,1,...\}$ on the vertex
set $V(G)$ with the transition probability matrix
\begin{equation}
\label{eq:transmat}
P=D^{-1}A,
\end{equation}
whose elements are
$$
p_{ij}=P[X_{t+1}=j|X_t=i]=\left\{\begin{array}{ll}
1/d_i, & \mbox{if} \ i \sim j,\\
0, & \mbox{otherwise}.
\end{array}\right.
$$
SRW is a time-reversible process with the stationary distribution
\begin{equation}
\label{eq:SRWstat}
\pi_i = \frac{d_i}{2m}, \quad i=1,...,n.
\end{equation}
Since SRW is time-reversible, its transition matrix is similar to a symmetric
matrix and hence the eigenvalues of the transition matrix are real, semi-simple and
can be indexed as follows:
$$
1=\lambda_1 \ge \lambda_2 \ge ... \ge \lambda_n \ge -1.
$$
Denote $\lambda_*$ the maximum modulus eigenvalue of $P$ different from $+1$ and $-1$
and call $\gamma(P)=1-|\lambda_*|$ the spectral gap.

The relaxation time $t_{rel}$ is then defined by
\begin{equation}
\label{eq:reltime}
t_{rel} = \frac{1}{\gamma(P)}.
\end{equation}
One interpretation of the relaxation time is as follows \cite{MTbook}:
if $t \ge t_{rel}$, then the standard deviation of $P^tf$ is bounded
by $1/e$ times the standard deviation of $f$. Also, the following
inequality (see, e.g., \cite{AF02,C97,MTbook}) indicates a strong relation between
the relaxation time and mixing time
$$
(\log(1/\varepsilon)+\log(1/2))(t_{rel}-1) \le t_{mix}(\varepsilon)
\le (\log(1/\varepsilon)+\log(1/\pi_{min})) t_{rel}.
$$
In particular, the above inequality suggests that for a finite Markov chain
and for small enough $\varepsilon$, the $\varepsilon$-mixing time is very close to
$\log(1/\varepsilon)t_{rel}$.

Now let us define the random walk with jumps (RWJ). It is a random walk
on a combination of the original graph and the complete graph weighted
by a scaled parameter $\alpha/n$ \cite{ART10}. Specifically, let us modify the adjacency matrix
in the following way
\begin{equation}
\label{eq:modifadj}
A(\alpha) = A + \frac{\alpha}{n} \one \one^T.
\end{equation}
Note that the new degree matrix is given by $D(\alpha)=D+\alpha I$, where $I$ is the identity
matrix. Then, the random walk on the modified graph is described by the following
transition probability matrix
\begin{equation}
\label{eq:modifadj}
P(\alpha) = D^{-1}(\alpha) A(\alpha),
\end{equation}
with elements
$$
p_{ij}(\alpha)=\left\{\begin{array}{ll}
\frac{1 + \alpha /n}{d_i +  \alpha}, & \mbox{if } i\sim j, \\
\frac{\alpha/n}{d_i +  \alpha}, & \mbox{otherwise},
\end{array}\right.
$$
Since RWJ is again a random walk on a weighted undirected graph, it is time-reversible Markov chain
with semi-simple eigenvalues. The stationary distribution of RWJ also has a simple form
\begin{equation}
\label{eq:RWJstat}
\pi_i(\alpha) = \frac{d_i+\alpha}{2m + \alpha n}, \quad i=1,...,n.
\end{equation}
The modified transition matrix $P(\alpha)$ can be rewritten as follows:
\begin{equation}
\label{eq:ndPR}
P(\alpha) = (D+\alpha I)^{-1}D \, P + (D+\alpha I)^{-1} \alpha I \, \one \left(\frac{1}{n}\one^T\right).
\end{equation}
We note that if the graph is regular, i.e., $D=d I$, the above expression reduces to
$$
P(\alpha) = \frac{d}{d+\alpha} P + \frac{\alpha}{d+\alpha}  \one \left(\frac{1}{n}\one^T\right),
$$
which is the transition matrix for PageRank (PR) \cite{PR} with the damping factor $d/(d+\alpha)$.
Thus, in the case of a regular graph RWJ coincides with PR. However, when the graph
has inhomogeneous degrees, these two concepts are different. From the equivalence of RWJ
to PR on the regular graph, we can immediately conclude that the relaxation time is
monotonously decreasing with $\alpha$ when the graph is regular. When the graph is irregular,
the situation becomes much more complex.

We also would like to note that RWJ can be viewed
as a node-dependent PageRank \cite{AHS14}, where the restart probability at each node is given by
$\alpha/(d_i+\alpha)$. Thus, in contrast to PR, RWJ restarts with smaller probabilities from
higher degree nodes.

\section{Application of Dobrushin coefficient}
\label{sec:dobr}

The Dobrushin ergodic coefficient (see, e.g., \cite{B99,IS11,S06}) can be used to obtain a lower bound
on the spectral gap $\gamma(P)$ of a Markov chain. The Dobrushin ergodic coefficient
is given by
\begin{equation}
\label{eq:dobr}
\delta (P) = \frac{1}{2} \max_{i,j \in V} \sum_{k \in V} |p_{ik}-p_{jk}|,
\end{equation}
or, equivalently,
\begin{equation}
\label{eq:dobr1}
\delta (P) = 1 - \min_{i,j \in V} \sum_{k \in V} p_{ik} \wedge p_{jk}.
\end{equation}
In the case of RWJ, we can obtain a simple upper bound on $\delta (P)$ by taking the smallest
element in the transition probability matrix. Namely, we obtain
\begin{equation}
\label{eq:transprob}
\delta (P) \le 1 - \min_{i,j \in V} \sum_{k \in V} \frac{\alpha/n}{d_{max}+\alpha}
= 1 - \frac{\alpha}{d_{max}+\alpha},
\end{equation}
where $d_{max}$ is the maximal degree in the graph. Since $\gamma(P) \ge 1-\delta(P)$,
we also obtain a lower bound for the spectral gap
\begin{equation}
\label{eq:dobr_bound}
\gamma(P(\alpha)) \ge \frac{\alpha}{d_{max}+\alpha}.
\end{equation}
And since $\alpha/(d_{max}+\alpha) \to 1$ as $\alpha \to \infty$, we have

\begin{pro}
For any undirected graph $G$, there always exists $\bar{\alpha}=\bar{\alpha}(G)$ such that
for all $\alpha > \bar{\alpha}$, we have $\gamma(P(\alpha)) > \gamma(P)$.
\end{pro}

For the regular graphs with degree $d$, the bound $\alpha/(d+\alpha)$ is quite tight.
In fact, for the regular graphs, as was noted at the end of the previous section
$\frac{\alpha}{d+\alpha}$ corresponds to the restart probability and
the exact value of the spectral gap is $\frac{d}{d+\alpha}\lambda_*(P)$ \cite{HK03,LM06}.
However, we have observed that for irregular graphs the bound (\ref{eq:dobr_bound}) can be very loose.

\section{Conditions in the case of small jump rate}
\label{sec:NDcond}

Let us analyse in this section the effect of small jump rate on the relaxation time.

Denote by $v_{*}(\alpha)$ the eigenvector corresponding to $\lambda_{*}(\alpha)$, that is
\begin{equation}
\label{eq:perteig}
P(\alpha) v_*(\alpha) = \lambda_{*}(\alpha) v_*(\alpha).
\end{equation}
For brevity, we shall write $v_*=v_{*}(0)$ and $\lambda_*=\lambda_{*}(0)$.
We also need the following preliminary result
\begin{lem}
The eigenelements $v_*(\alpha)$ and $\lambda_*(\alpha)$ are analytic functions with
respect to $\alpha$.
\end{lem}
{\bf Proof:} It is known from \cite[Chapter 2]{Kato} (see also \cite{AFH13,B85}) that
if the eigenvalues of the perturbed matrix are semi-simple, then the eigenvalues as well as
eigenvectors can be expanded as power series with positive integer powers. Since in our case,
RWJ is time-reversible, its eigenvalues are semi-simple and the statement of the lemma follows. \hfill $\Box$

\bigskip

Next, we are in a position to provide necessary and sufficient conditions for the
improvement of the relaxation time for small jump rates.

\begin{theorem}
For sufficiently small $\alpha$, in the case $\lambda_{*} < 0$,
the spectral gap increases, and equivalently, the relaxation time decreases with respect to $\alpha$.

If $\lambda_{*} >0$, for sufficiently small $\alpha$ the necessary and sufficient condition for the
decrease of the relaxation time is
\begin{equation}
\label{eq:NandS}
    \frac{1}{n}(\one^T v_{*})^2  < \lambda_{*} v_{*}^Tv_{*}.
\end{equation}

In addition, we have the following asymptotics:
$$
\lambda_{*}(\alpha)  = \lambda_{*}(0) + \alpha \frac{\frac{1}{n}(v^{(0)}_{*})^T \one\one^T v^{(0)}_{*} - \lambda^{(0)}_{*}(v^{(0)}_{*})^Tv^{(0)}_{*}}{{(v^{(0)}_{*})^T Dv^{(0)}_{*}}} + O(\alpha^2).
$$
\end{theorem}
{\bf Proof:}
According to Lemma~1, the eigenelements $\lambda_{*}(\alpha)$ and $v_*(\alpha)$ can be expanded
as power series
\begin{equation}
\label{eq:Taylorlambda}
\lambda_{*}(\alpha)  = \lambda^{(0)}_{*} + \alpha \lambda^{(1)}_{*} + O(\alpha^2),
\end{equation}
\begin{equation}
\label{eq:Taylorv}
v_{*}(\alpha) = v^{(0)}_{*} + \alpha v^{(1)}_{*} + O(\alpha^2),
\end{equation}
for sufficiently small $\alpha$. If we set $\alpha=0$, we obtain
$\lambda^{(0)}_{*}=\lambda_{*}(0)=\lambda_{*}$ and
$v^{(0)}_{*}=v_{*}(0)=v_{*}$ such that
$$
D^{-1}A v_{*} = \lambda_{*} v_{*},
$$
or equivalently,
\begin{equation}
\label{eq:eig0}
A v_{*} = \lambda_{*} D v_{*}.
\end{equation}
It is also convenient to rewrite (\ref{eq:perteig}) as a generalized eigenvalue problem
\begin{equation}
\label{eq:genperteig}
(A + \frac{\alpha}{n} \one\one^T) v_*(\alpha) = \lambda_{*}(\alpha) (D + \alpha I)  v_*(\alpha).
\end{equation}
Substituting the power series (\ref{eq:Taylorlambda}) and (\ref{eq:Taylorv}) into
equation (\ref{eq:genperteig}) and equating coefficients in $\alpha$-terms, i.e.,
$$
(A + \frac{\alpha}{n} \one\one^T) (v^{(0)}_{*} + \alpha v^{(1)}_{*} + O(\alpha^2)) =
$$
$$
(\lambda^{(0)}_{*} + \alpha \lambda^{(1)}_{*} + O(\alpha^2)) (D + \alpha I) (v^{(0)}_{*} + \alpha v^{(1)}_{*} + O(\alpha^2)),
$$
yields
\begin{align}
\frac{1}{n}\one\one^Tv^{(0)}_{*} + Av^{(1)}_{*} = \lambda^{(0)}_{*}Iv^{(0)}_{*} + \lambda^{(1)}_{*}Dv^{(0)}_{*} + \lambda^{(0)}_{*}Dv^{(1)}_{*}.
\end{align}
Now let us multiply the above equation by $v^T_{*}$ from the left to obtain
\begin{equation}
\label{eq:longone}
\frac{1}{n}(v^{(0)}_{*})^T \one \one^T v^{(0)}_{*} + (v^{(0)}_{*})^T Av^{(1)}_{*} =
\end{equation}
$$
(v^{(0)}_{*})^T\lambda^{(0)}_{*}v^{(0)}_{*} + (v^{(0)}_{*})^T\lambda^{(1)}_{*}Dv^{(0)}_{*} +  (v^{(0)}_{*})^T\lambda^{(0)}_{*}Dv^{(1)}_{*}.
$$
Due to symmetry, the equation (\ref{eq:eig0}) can be rewritten as
$$
v_{*}^T A = \lambda_{*} v_{*}^T D,
$$
and hence
$$
(v^{(0)}_{*})^T Av^{(1)}_{*} = (v^{(0)}_{*})^T\lambda^{(0)}_{*}Dv^{(1)}_{*},
$$
which simplifies (\ref{eq:longone}) to
$$
\frac{(v^{(0)}_{*})^T \one\one^Tv^{(0)}_{*}}{n}  = (v^{(0)}_{*})^T\lambda^{(0)}_{*}v^{(0)}_{*} + (v^{(0)}_{*})^T\lambda^{(1)}_{*}Dv^{(0)}_{*}.
$$
Thus, $\lambda^{(1)}_{*}$ can be expressed as
\begin{equation}
\label{eq:lambdaexpressor}
\lambda^{(1)}_{*} = \frac{\frac{1}{n}(v^{(0)}_{*})^T \one\one^T v^{(0)}_{*} - \lambda^{(0)}_{*}(v^{(0)}_{*})^Tv^{(0)}_{*}}{{(v^{(0)}_{*})^T Dv^{(0)}_{*}}}.
\end{equation}
Clearly, the denominator in (\ref{eq:lambdaexpressor}) is always positive. Now, consider two cases:

\begin{case} $\lambda_{*} < 0$
\end{case}
In this case the numerator in (\ref{eq:lambdaexpressor}) is always positive. Then $\lambda^{(1)}_{*}$ is also positive.
By expansion (\ref{eq:Taylorlambda}), when $\alpha$ is sufficiently small, the absolute value of $\lambda_*$
is decreasing with respect to $\alpha$.

\begin{case} $\lambda_{*} > 0$
\end{case}
Again, by expansion (\ref{eq:Taylorlambda}), for sufficiently small $\alpha$ the value of $\lambda_{*}(\alpha)$
decreases in $\alpha$ if and only if  $\lambda^{(1)}_{*}$ is negative, i.e., when the numerator is negative.
This is precisely what was stated in the theorem's condition. \hfill $\Box$

\bigskip

\section{Sufficient conditions with easier interpretation}
\label{sec:nodedegcond}

Even though the condition in Theorem~1 is necessary and sufficient, it is not easy to use and does not
have an easy intuitive interpretation. Next, we shall derive a series of sufficient conditions with easier interpretation
and verification. Towards this goal, let us transform the condition in Theorem~1 to
an equivalent form using the combinatorial Laplacian. Specifically, we shall use the combinatorial
Laplacian of the complete graph:
$$
L_K = nI-\one \one^T.
$$

\begin{lem}
The condition $(\ref{eq:NandS})$ is equivalent to
$$
1 - \lambda_* < \frac{v_*^T L_K v_*}{n v_*^T v_*},
$$
and if $\lambda_* > 0$, the condition $(\ref{eq:NandS})$ is equivalent to
\begin{equation}
\label{eq:NandD1}
\gamma(P)  < \frac{\sum_{i, j}(v_{*i} - v_{*j})^2}{n \sum_{i} v^2_{*i}}.
\end{equation}
\end{lem}
{\bf Proof:}
Using the definition of the Laplacian, we can write
$$
\frac{1}{n} (\one^T v_*)^2 = \frac{1}{n} v_*^T \one \one^T v_* = \frac{1}{n} v_*^T (nI - L_K) v_*
= v_*^Tv_* - \frac{1}{n} v_*^T L_K v_*.
$$
Thus, we can rewrite condition (\ref{eq:NandS}) in a new form:
$$
1 - \lambda_* < \frac{v_*^T L_K v_*}{n v_*^T v_*}.
$$
The other equivalent form follows immediately from the definitions of the spectral gap and $L_K$.
\hfill $\Box$

\bigskip

Next we provide a couple of sufficient conditions with easy interpretation in the context of complex networks.

\begin{cor}
If
\begin{equation}
\gamma(P) < \frac{1}{n},
\end{equation}
then, for sufficiently small $\alpha>0$, the spectral gap of $P(\alpha)$ is larger than the spectral gap of $P$.
\end{cor}
{\bf Proof:}
Clearly, $\sum \limits_{i,j} (v_{*i}(\alpha)-v_{*j}(\alpha))^2 > \sum \limits_{i} v_{*i}(\alpha)^2$.
Because $v_{*}(\alpha) \perp \pi(\alpha)$ and so there are both positive and negative numbers among $v_{*i}(\alpha)$.
Each number $v_{*i}(\alpha)$ has at least one number of opposite sign $v_{*j}(\alpha)$ and such that
$(v_{*i}(\alpha) - v_{*j}(\alpha))^2 > v_{*i}(\alpha)^2$. \hfill $\Box$

\bigskip

The above corollary has the following simple and useful interpretation: If $P$ has a sufficiently small gap,
the addition of jumps with small rate always improves the relaxation time. It is known that many complex networks have very
small spectral gap and thus this corollary gives an explanation why the relaxation time typically improves
with the addition of jumps in complex networks. The next corollary refines the above argument.

\bigskip

\begin{cor}
Denote the proportion of negative and positive $v_{*i}(\alpha)$ as $\mu$ : $1 - \mu$.
Then, for sufficiently small $\alpha$, if
\begin{equation}
\gamma(P) < \min(\mu, 1 - \mu),
\end{equation}
the relaxation time decreases with respect to $\alpha$.
\end{cor}
{\bf Proof:}
For each positive number $v_{*i}(\alpha)$ has $\mu n$ numbers of opposite sign $v_{*j}(\alpha)$ such that $(v_{*i}(\alpha) - v_{*j}(\alpha))^2 > v_{*i}(\alpha)^2$. Analogously, each negative number $v_{*i}(\alpha)$ has $(1 - \mu)n$ numbers of opposite sign $v_{*j}(\alpha)$ and for them $(v_{*i}(\alpha) - v_{*j}(\alpha))^2 > v_{*i}(\alpha)^2$. So, $\sum \limits_{i,j} (v_{*i}(\alpha)-v_{*j}(\alpha))^2 \ge n \min(\mu, 1 - \mu) \sum \limits_{i} v_{*i}(\alpha)^2$. \hfill $\Box$

\bigskip

As a result of Corollary~2, the closer $\mu$ is to $1/2$, the better. Often complex networks have clustering
structure. The eigenvector $v_{*}$ can be interpreted as a variant of the Fiedler vector. Thus, if a complex network
can be divided into two clusters of similar sizes, the value of $\mu$ will be close to 1/2 or at least far from zero.
In such a case, the spectral gap of the original transition matrix $P$ does not need to be small for the addition
of jumps to improve the relaxation time. The above statement is similar in spirit to the condition given in \cite{ART10}.

\bigskip

Let us now provide a sufficient condition in terms of nodes' degrees.

\bigskip

\begin{theorem}
Let the following condition hold:
\begin{equation} \label{degrees_averages}
\gamma(P) < 4 \frac{(d_1+d_2+\dots+d_n)^2}{n(d_1^2+d_2^2+\dots+d_n^2)}.
\end{equation}
Then, if $\alpha > 0$ is sufficiently small,  $\gamma(P(\alpha)) > \gamma(P)$.
\end{theorem}
{\bf Proof:}
We replace the condition (\ref{eq:NandD1}) with the following more stringent condition:
\begin{equation}
\label{rewritten_main_equation}
\gamma(P) < \min_{\vec{f} \bot \vec{\pi}} \frac{\sum_{i \sim j} (f_i-f_j)^2}{n\sum{f_i^2}}.
\end{equation}
Let us find the minimum of the RHS of the above expression over all vectors $\vec{f}$, orthogonal to $\pi(\alpha)$,
i.e., for vectors satisfying $\sum {f_i}{d_i} = 0$.
So it will be enough for $\gamma(P)$ to be less than the RHS of the condition (\ref{rewritten_main_equation}).

Note that the RHS of (\ref{rewritten_main_equation}) is homogeneous in $\vec{f}$, so without loss
of generality we can set $\sum {f_i}^2 = 1$.

We shall use the method of Lagrange multipliers for the objective function $\sum_{i, j}(f_i - f_j)^2$ and the following constraints:
\begin{equation} \label{squares_sum_one}
\sum_{i=1}^n f_i^2 = 1,
\end{equation}
\begin{equation}\label{weighted_sum_zero}
\sum_{i=1}^n {f_i}{d_i} = 0.
\end{equation}
The expression $\sum_{i, j}(f_i - f_j)^2$ taking into account the latter constraint can be rewritten as $2n - 2(f_1+f_2+\cdots+f_n)^2$.
Therefore we are looking for the maximum of $(f_1+f_2+\cdots+f_n)^2$, which will give us the minimum of $\sum_{i, j}(f_i - f_j)^2$.
Let us find the maximal absolute values of the function $(f_1+f_2+\cdots+f_n)$ with respect to the above constraints. Obviously, these extreme values exist (sums cannot be more than $n$ and less than $-n$). Furthermore, their absolute values coincide, because domains of these two sums is symmetric with the center in the origin. So the set of values of $(f_1+f_2+\cdots+f_n)$ is centrally symmetric and the maximal absolute values of $(f_1+f_2+\cdots+f_n)$ and $-(f_1+f_2+\cdots+f_n)$ are the same.

Looking at the geometry of our maximization problem, we can see that the constraints give us a sphere cut by the hyperplane $\sum {f_i}{d_i} = 0$. We are looking for the touching points of this set and some hyperplane $(f_1+f_2+\cdots+f_n)=c$. Clearly, there will be exactly two centrally-symmetric points of touching, giving us the maxima of the absolute values (the domain of this sum is symmetric with the center
of the symmetry in the origin, so the set of values of $(f_1+f_2+\cdots+f_n)$ is centrally symmetric).

Now, consider the Lagrange function:
\begin{equation}
L = (f_1+f_2+\cdots+f_n)  - \lambda_1 \Big(\sum_{i} f_i^2 - 1\Big) - \lambda_2\sum {f_i}{d_i}.
\end{equation}

Taking partial derivatives with respect to $f_i$ for $i=1,...,n$, we obtain
\begin{equation}
0 = 1  - 2\lambda_1 f_i - \lambda_2 d_i.
\end{equation}
Firstly, multiply all these equations with $d_i$, respectively, and sum all of them (here we also use equation (\ref{weighted_sum_zero})). The result is
\begin{equation}
0 = 2m  - \lambda_2 \sum_{i} d_i^2.
\end{equation}
Thus, $\lambda_2 = \frac{2m}{\sum_{i} d_i^2}$, where $m$ is the total number of edges.

Secondly, multiply all equations with $f_i$, respectively, and sum all of them (here we additionally use equations (\ref{squares_sum_one}) and (\ref{weighted_sum_zero})). The result is
\begin{equation}
0 = S  - 2\lambda_1,
\end{equation}
where $S = (f_1+f_2+\cdots+f_n)$. Thus, $\lambda_1 = S/2$.

And finally, just sum all equations (here we also use equations (\ref{squares_sum_one}) and (\ref{weighted_sum_zero})) to get
\begin{equation}
0 = n  - 2\lambda_1 S - \lambda_2 \big(2m\big).
\end{equation}
If we know $\lambda_1$ and $\lambda_2$, we also know $S^2 = n - \frac{8m^2}{\sum_{i} d_i^2}$.
Substituting this into our initial expression that we want to minimize, we obtain
\begin{equation}
\min_{\vec{f} \bot \vec{\pi}} \frac{\frac{1}{n} \sum_{i, j} (f_i-f_j)^2 }{\sum{f_i^2}} =  \frac{1}{n} \Big( 2n -2\big(n -  \frac{8m^2}{\sum_{i} d_i^2}\big) \Big) = 4\frac{(d_1+d_2+\dots+d_n)^2}{n(d_1^2+d_2^2+\dots+d_n^2)}.
\end{equation}

Thus, if inequality (\ref{rewritten_main_equation}) holds, the condition of Theorem~1 also holds,
and hence $\gamma(P(\alpha)) > \gamma(P)$ for sufficiently small $\alpha$. \hfill $\Box$

\bigskip

We can also rewrite the condition of the above theorem as follows:

\bigskip

\begin{cor}
Let the following condition hold:
\begin{equation} \label{degrees_averages}
\gamma(P) < 4 \frac{(\bar{d})^2}{\overline{d^2}},
\end{equation}
where
$$
\bar{d} = \frac{d_1+d_2+\dots+d_n}{n},
$$
and
$$
\overline{d^2} = \frac{d_1^2+d_2^2+\dots+d_n^2}{n}.
$$
Then, if $\alpha > 0$ is sufficiently small,  $\gamma(P(\alpha)) > \gamma(P)$.
\end{cor}

The quantity $(\bar{d})^2/\overline{d^2}$ is the reciprocal of the efficiency or of the squared coefficient of variation.
It is also sometimes referred to as the signal-to-noise ratio. From Corollary~3 we can see that the more ``irregular''
the degree sequence is, the more stringent is the condition on the spectral gap.

\bigskip

The following simple sufficient condition using only $\bar{d}$ and $d_{max}$ also holds.

\bigskip

\begin{cor}
For sufficiently small $\alpha$, $\gamma(P(\alpha)) > \gamma(P)$ is true for all graphs such that
\begin{equation}
\gamma(P) < 4\bar{d}/d_{max}.
\end{equation}
\end{cor}
{\bf Proof:}
The RHS of the inequality (\ref{degrees_averages}) is greater than $\bar{d}/d_{max}$. \hfill $\Box$

\section{Examples}
\label{sec:num}

Let us first demonstrate in this section that there exist weighted graphs for which the introduction
of jumps actually increases the relaxation time. Towards this end, we consider a weighted graph of size 2
with the following adjacency matrix:
$$
A=\left[\begin{array}{cc}
a_{11} & a_{12}\\
a_{12} & a_{22}
\end{array}\right].
$$
The characteristic equation for the generalized eigenvalue problem $Av=\lambda Dv$ takes the form
$$
\det(A-\lambda D)=
\det \left[\begin{array}{cc}
a_{11}-\lambda(a_{11}+a_{12}) & a_{12}\\
a_{12} & a_{22}-\lambda(a_{22}+a_{12})
\end{array}\right]=0.
$$
It has two solutions
$$
\lambda_1=1 \quad \mbox{and} \quad \lambda_*=\frac{\det(A)}{(a_{11}+a_{12})(a_{22}+a_{12})}.
$$
The eigenvector corresponding to $\lambda_*$ is given by
$$
v_*=
\left[\begin{array}{c}
1\\
-\frac{a_{11}+a_{12}}{a_{22}+a_{12}}
\end{array}\right]C.
$$
Let us calculate the numerator in the expression (\ref{eq:lambdaexpressor})
$$
\frac{1}{2} (v_*^T \one)^2 - \lambda_* v_*^T v_* =
$$
$$
\frac{(a_{22}-a_{11})^2(a_{11}+a_{12})(a_{22}+a_{12})-2\det(A)[(a_{11}+a_{12})^2+(a_{22}+a_{12})^2]}
{2(a_{11}+a_{12})(a_{22}+a_{12})^3}.
$$
Now if we choose the weights such that $\det(A)$ is close to zero or just zero and $a_{11} \neq a_{22}$,
the spectral gap will decrease, and consequently the relaxation time will increase when $\alpha$ increases
from zero. Thus, in general, the conditions of Theorem~1 need to be checked in the case of weighted graphs
if one wants to be sure that the relaxation time decreases with the increase of the jump rate.

We also tried to construct an example of unweighted graph when the relaxation time increases with the
increase in the jump rate. Somewhat surprisingly, it appears to be hard to construct such an example.
In fact, we have checked all non-isomorphic graphs of sizes up to $n \le 9$ from the Brendan McKay's collection \cite{BMKCol}
and could not find a single example when the relaxation time increases with the introduction of small
jump rate. The code of this verification can be found at \cite{Code}.

We have also checked various random graph models like Erd\H{o}s-R\'enyi graph or non-homogeneous
Stochastic Block Model and could not find an example of an unweighted graph for which the relaxation time
increases with the introduction of jumps.

\section{Conclusion and future research}
\label{sec:conc}

We have analysed the spectral gap, or equivalently the relaxation time, in the random walk with jumps.
We have obtained a necessary and sufficient condition for the decrease of the relaxation time when
the jump rate increases from zero. These conditions unfortunately do not have easy interpretation.
Therefore, we have proceeded with the derivation of several sufficient conditions with easy interpretation.
Some of these sufficient conditions can also be easily verified.
The derived sufficient conditions show that in most complex networks the relaxation time should
decrease with the introduction of jumps.
We have also demonstrated that there exist
weighted graphs for which the relaxation time increases with the introduction of jumps. On the other
hand, we could not find such an example in the case of unweighted graphs. At the moment, we tend to
conjecture that the introduction of jumps always improves the relaxation time in the case of unweighted
graphs.

\section*{Acknowledgements}
The work was partially supported by the joint laboratory Nokia Bell Labs - Inria.
This is an author version of the paper accepted to WAW 2018.

\end{document}